\documentclass{article}

\usepackage{delarray,verbatim,enumerate,a4}
\usepackage{amsmath,amsthm,amstext,amsbsy,amssymb,amsfonts,amscd}
\usepackage[all]{xy}
\usepackage[frenchb]{babel}

\usepackage{aeguill}

\newtheorem{Th}{Théorème}[]
\newtheorem{Lem}[Th]{Lemme}
\newtheorem{Prop}[Th]{Proposition}
\newtheorem{Cor}[Th]{Corollaire}
\newtheorem{Conj}[Th]{Conjecture}
\newtheorem{Sco}[Th]{Scolie}
\newtheorem{Def} [Th]{Définition}

\def\Preuve{\smallskip\noindent {\it Preuve.~}}
\def\PreuveTh{\smallskip\noindent {\it Preuve du Théorème.~}}
\def\Remarque{\smallskip\noindent {\it Remarque.~}}
\def\Remarques{\smallskip\noindent {\it Remarques.~}}
\def\Exemple{\smallskip\noindent {\it Exemple.~}}

\font\teneufm=eufm10
\font\seveneufm=eufm7
\font\fiveeufm=eufm5
\newfam\gothfam
\textfont\gothfam=\teneufm \scriptfont\gothfam=\seveneufm
\scriptscriptfont\gothfam=\fiveeufm
\def\goth{\fam\gothfam}

\def\Fl{\mathbb{F}_\ell}		\def\Q{\mathbb Q}		
\def\Z{\mathbb Z}			\def\Zl{{\mathbb{Z}_\ell}} 
\def\Tl{\mathbb{T}_\ell} 		\def\N{\mathbb N}
\def\Ql{{\mathbb{Q}_\ell}} 	\def\lT{\ov{\mathbb T}_\ell}

\def\wi{\widetilde}			\def\ov{\overline}

 		\def\U{\mathcal  U}	\def\E{\mathcal  E} 
\def\J{\mathcal  J}  	\def\C{\mathcal  C}	\def\R{\mathcal  R}
\def\Dl{\mathcal  D\ell} 	\def\Pl{\mathcal  P\ell}  	\def\Cl{\mathcal  C \ell}

 	\def\W{W\!K}	   	  \def\T{\mathcal  T}

		\def\p{{\goth p}}		\def\l{{\goth l}}
\def\d{{\goth d}}		\def\P{{\goth  P}}

\def\div{\operatorname{div}}	\def\deg{\operatorname{deg}}
\def\Gal{\operatorname{Gal}}	\def\Rad{\operatorname{Rad}}
\def\Ker{\operatorname{Ker}}	
\def\dim{\operatorname{dim}}	\def\rg{\operatorname{rg}}
\def\Log{\operatorname{Log}}	\def\Hom{\operatorname{Hom}}

\begin{document}


\title{\Large \bf Approche logarithmique des noyaux étales sauvages 
des corps de nombres\footnote{J. Number Th. {\bf 120}  (2006), 72--91.}}

\author{Jean-François {\sc Jaulent}  {\small \&} Alexis {\sc Michel}}

\date{}
\maketitle

{\small
\noindent{\bf Abstract.}
 We study the $\ell$-part of the the wild {\it étale} kernels $\W_{2i}(F)$
of an arbitary number field $F$ for a given prime $\ell$  in connection 
with the logarithmic $\ell$-class groups $\wi\Cl_F$. From the logarithmic
arithmetic we deduce rank formulas, periodicity and reflection theorems,
triviality characterizations and various consequences.

\medskip\noindent
{\bf Résumé.} Nous étudions la $\ell$-partie des noyaux étales sauvage 
$\W_{2i}(F)$ d'un corps de nombres arbitaire $F$ en liaison avec 
l'arithmétique des classes logarithmiques. Nous en déduisons 
notamment des formules de rang, des résultats de périodicité et de 
réflexion, des caractérisations de la trivialité, ainsi que diverses 
conséquences.}\bigskip

\bigskip

\noindent{\large \bf  Introduction.}

\bigskip

Introduits par P. Schneider \cite{Sc} comme noyaux de localisation pour certains groupes de cohomologie $\ell$-adique attachés à un corps de nombres $F$, les $\ell$-noyaux étales sauvages $\W_{2i}(F)$ ont été interprétés en termes de $K$-théorie supérieure suite en particulier aux travaux de C. Soulé \cite{So}, W. Dwyer \& E. Friedlander \cite{DF}, G. Banaszak \cite{Ba1, Ba2,Ba3}, et ainsi nommés par T. Nguyen Quang Do \cite{N0} car ils coïncident pour $2i = 2$ avec les noyaux sauvages usuels de la $K$-théorie, en vertu de la conjecture de Quillen-Lichtenbaum établie dans ce cas par Tate.\smallskip

Nous sommes plus particulièrement intéressés ici par certaines de leurs propriétés arithmétiques déjà étudiés par plusieurs auteurs~; notamment par M. Kolster \cite{K3}, en liaison avec la conjecture de Leopoldt ~; par 
J.-F. Jaulent et F. Soriano \cite{JS2} pour $2i=2$, sous celle de Gross (cf. \cite{FG,J6}), en lien avec le $\ell$-groupe des classes logarithmiques.\smallskip

L'objet du présent travail est ainsi de généraliser les résultats de \cite{JS2} à {\em n'importe quel} corps de nombres, pour tout  $i$ dans $\Z$ et {\em indépendamment de toute conjecture}.\smallskip

Le principe de notre méthode est simple~: Partant de la description 
tordue {\it à la Tate} de ces divers groupes en haut de la tour cyclotomique, 
nous en déduisons par codescente, sous l'hypothèse essentielle de procyclicité\footnote{Hypothèse qui est naturellement vérifiée pour tous les premiers $\ell$ impairs, mais non pour $\ell=2$.} de la tour cyclotomique $F[\zeta_{\ell^\infty}]/F$,
des relations précises entre les $\ell$-noyaux  $\W_{2i}(F)$ et le $\ell
$-groupe des classes logarithmiques du corps $F[\zeta_{2\ell}]$, lesquelles 
permettent en retour de transporter aux noyaux étales sauvages les 
informations sur ce groupe données par la Théorie logarithmique des 
genres et, plus généralement, l'arithmétique des classes et des unités 
logarithmiques développée {\it pour elle-même} dans \cite{J6}.\smallskip

Outre les perspectives algorithmiques qu'ouvre notre approche (cf. \cite{DS,Dy}), il nous
parait important de souligner trois points~: D'abord la simplicité des 
résultats obtenus\footnote{Le Corollaire 5 en est un exemple parfait.}. 
Ensuite les implications heuristiques de cette description, dont la 
discussion finale en appendice est une assez belle illustration. Enfin l'importance de l'hypothèse de procyclicité pour $\ell=2$, en l'absence de laquelle le {\em bon} approximant des 2-noyaux étales n'est plus le 2-groupe des classes logarithmiques qui apparait dans \cite{J5}, mais le 2-groupe des classes positives $\Cl^{\,pos}_F$ introduit à cet effet dans \cite{JS3}, en cohérence avec les résultats de \cite{H1}. Ce dernier point est à mettre en relation avec le fait que le 2-noyau sauvage $\W_{2i}(F)$ ne coïncide pas dans ce cas avec le sous-groupe des éléments de hauteur infinie de la 2-partie du groupe $K_{2i}(F)$ (cf. \cite{Ta,H2,H3,HR,JS4}).
\bigskip

{\small 
\noindent{\bf \em Remerciements.} Le présent article est la rédaction d'une conférence donnée à Metz en juin 2002 lors d'un colloque sur la $K$-théorie des corps de nombres. Les auteurs remercient tout particulièrement le rapporteur pour les améliorations apportées à leur texte.}
\smallskip


\bigskip

\noindent{\large \bf 1. Classes logarithmiques et noyaux étales sauvages}

\bigskip


Nous nous plaçons tout au long de cet article dans la situation générale suivante~: le nombre premier $\ell$  étant fixé, nous désignons par $F$ un corps de nombres arbitraire, à cette seule réserve\footnote{Cette condition n'est évidemment restrictive que pour $\ell = 2$.} que nous demandons que l'extension cyclotomique $F[\zeta_{\ell^\infty}]/F$, engendrée par les racines d'ordre $\ell$-primaire de l'unité, soit {\it procyclique}.


\bigskip

\noindent{\bf 1.a. Le contexte galoisien semi-simple}

\medskip


Le corps $F$ étant fixé, nous notons $L$ une extension abélienne de $F$, de 
groupe de Galois $\Delta\, =\, \Gal (L/F)$, de degré $d$ étranger à $\ell$, 
et contenant une racine pri\-mitive $\ell$-ième de l'unité $\zeta$. Sous
l'hypothèse $\ell \nmid d$, l'algèbre résiduelle $\Fl [\Delta ]$ est ainsi
une algèbre semi-simple, produit direct d'extensions $F_\varphi $ de 
$\Fl$~; l'algèbre $\ell$-adique $\Zl [\Delta ]$ une algèbre semi-locale, 
produit direct d'extensions non ramifiées $Z_\varphi $ de $\Zl$~;
et les idempotents primitifs $\bar e_\varphi$  (respectivement $e_\varphi $)
correspondant à leurs décompositions respectives~:
$$
\Fl [\Delta ]\, =\, \oplus _\varphi \ \Fl [\Delta ]\bar e_\varphi\,
 =\, \oplus _\varphi \ F_\varphi \ \ \ \ \ \& \ \ \ \ \ 
\Zl [\Delta ] \, =\, \oplus _\varphi \ \Zl [\Delta ] e_\varphi\, =
\, \oplus _\varphi \ Z_\varphi
$$
sont donnés à partir des caractères $\ell$-adiques irréductibles 
$\varphi $ de $\Delta $ par les formules classiques~:\smallskip

\centerline{$e_\varphi\, =\, \frac{1}{d}\sum_{\tau\in \Delta}
\varphi(\tau^{-1})\tau$, et leurs réductions respectives 
modulo  $\ell$.}\smallskip

\noindent Parmi les caractères de $\Delta $ figurent en particulier le
{\it caractére unité} 1, dont l'idempotent associé est donné à partir de 
la norme algébrique $\nu _\Delta \, =\,\sum_{\tau\in \Delta}\tau $ par 
$e_1\, =\, \frac{1}{d}\ \nu _{\Delta }$, et le {\it caractère cyclotomique} 
$\omega$, caractérisé\footnote{En particulier, on a l'égalité 
$\,\omega =1$ si et seulement si le corps $K$ contient $\zeta$.} par 
l'identité~:\smallskip

\centerline{$\zeta ^\tau\, =\, \zeta ^{\omega (\tau )}\quad \forall\, \tau\, \in
 \Delta$.}\smallskip

\noindent L'inverse\footnote{ De façon générale, il est commode de 
noter $\psi ^{-1}$ le caractère $\tau\mapsto \psi (\tau ^{-1})$}
 $\bar\omega =\omega ^{-1}$ de $\omega $, c'est à 
dire le caractère défini par $\bar\omega (\tau )=\omega (\tau ^{-1})$, 
est dit souvent $anticyclotomique$, et l'involution\smallskip

\centerline{$\psi\mapsto\psi ^* =\omega\psi ^{-1}$,}\smallskip

\noindent de l'algèbre $R_{\Z_{\ell} }(\Delta )$ des caractères $\ell$-adiques virtuels
de $\Delta $ est connue traditionnellement sous le nom d'{\it involution$ $ 
du miroir}~; on dit encore que $\psi ^* $ est le $reflet$  de $\psi$. Les
anneaux locaux $\Z_{\varphi^*}$, $\Z_{\varphi\omega^i}$ et $\Z_{\varphi^*
\omega^i}$, pour $i \in \mathbb N$, sont évidemment isomorphes à 
$\Z_\varphi$~; ils sont généralement notés $\Z_\phi$ dans ce qui suit.\medskip

Tous les $\Zl [\Delta ]$-modules que nous sommes amenés à considérer
plus loin pro\-viennent de $\Zl$-modules galoisiens fonctoriellement 
attachés à un corps de nombres. Il existe donc des applications 
naturelles de transition, ainsi dans $L/F$
la norme arithmétique $N_{L/F}$  (ou restriction) et le morphisme
d'extension $j_{L/F}$ (ou corestriction), qui vérifient les identités~:\smallskip

\centerline{$N_{L/F}\,\circ \, j_{L/F}=\ [L:F]$ \quad {\small \&} \quad  
$j_{L/F}\, \circ\, N_{L/F}=\nu _\Delta=\sum_{\tau\in \Delta}\tau$,}\smallskip

\noindent où $\nu_\Delta$ est la norme algébrique attachée au 
groupe $\Delta$. En particulier, l'ordre $|\Delta|=\ [L:F]$ de $\Delta$ étant par hypothèse 
inversible dans $\Zl$, il en résulte que le module $X_F$ associé au 
sous-corps $F$ de $L$ s'identifie canoniquement à l'image $X_L^{e_
\Delta}$ de celui attaché à $L$ par le projecteur associé à l'idempotent 
$e_\Delta=\frac{1}{|\Delta|}\ \nu_\Delta=\sum _{\varphi \in R_{\Zl}^
{\rm irr}(\Delta )} \ e_\varphi $~;
propriété que nous utilisons systématiquement dans ce qui suit.
\smallskip

\Exemple Pour illustrer immédiatement tout cela, regardons plus attentivement ce qui se passe dans la situation
non triviale la plus simple~: celle du miroir de  Scholz (cf. \cite{Gr}). 
Prenons donc $\ell =3 \ $, $F=\Q $ et considérons un corps biqua\-dratique
$L=\Q [\sqrt{-3},\sqrt d]$ (avec $d\in \mathbb N -\{ 0 ;1\}$ sans facteur 
carré). Le diagramme des sous-corps se présente alors comme suit~:
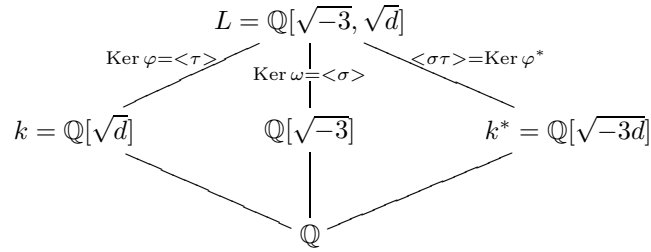
\smallskip

\begin{displaymath}
\xymatrix{{} & L = \Q [\sqrt{-3},\sqrt d] 
\ar@{-}[dl]_{\Ker\varphi = <\tau>} \ar@{-}[d]|{\Ker \omega=<\sigma>}
\ar@{-}[dr]^{\ <\sigma\tau> = \Ker \varphi^*}	\\ 
k=\Q [\sqrt{d}] \ar@{-}[dr] &\Q [\sqrt{-3}] \ar@{-}[d] & 
k^* =\Q [\sqrt{-3d}] \ar@{-}[dl]  \\ {}  & \Q}
\end{displaymath}
\smallskip

\noindent
Le groupe de Galois $\Delta  = \Gal (L/\Q )$ est le Vierergrüppe de  
Klein $V_4\, =\, C_2\times C_2$. Il possède  exactement quatre caractères 
3-adiques irréductibles~: le cara\-ctère unité 1, le caractère cyclotomique
$\omega $ (dont le noyau $\Ker\omega $ fixe précisément le sous-corps
cyclotomique $\Q [j]\, =\, \Q [\sqrt{-3}]$) et deux autres caractères de degré 1,
l'un réel $\varphi $ dont le noyau $\Ker \varphi $ fixe le sous-corps réel
$\Q [\sqrt d]$, l'autre imagi\-naire $\varphi ^*, =\, \omega \varphi
 \, =\, \bar \omega \varphi $ dont le noyau  $\Ker \varphi ^*$ fixe
le sous-corps quadratique imaginaire
$k^* =  \Q [\sqrt {-3d}]$. Si donc $\tau $ désigne la conjugaison
 complexe et $\sigma $ l'unique élément non trivial de $\Delta $ qui 
fixe le corps cyclotomique $\Q [j]$, les quatre idempotents primitifs de 
l'algèbre $\Z _3 [\Delta ]$
 sont respectivement~:
\begin{eqnarray*}
 e_1 &=& \frac{1}{4}(1+\tau )(1+\sigma ),  \qquad e_\omega  =
\frac{1}{4}(1-\tau )(1+\sigma ),\\
 e_\varphi & =& \frac{1}{4}(1+\tau )(1-\sigma ), \qquad
e_{\varphi ^*}\, =\, \frac{1}{4}(1-\tau )(1-\sigma )\ ;
\end{eqnarray*}
et les idempotents normiques associés aux trois sous-corps
quadratiques 
$\Q [j]$, \\ $k=\Q [\sqrt{d}]$ et $k^* =\Q [\sqrt{-3d}]$ 
sont ainsi~:\smallskip

\centerline{$\frac{1}{2}(1+\sigma )\, =\, e_1+e_{\omega} , \ \ \ \ \ \ 
\frac{1}{2}(1+\tau )\, =\, e_1+e_{\varphi} , \ \ \ \ \ \ 
\frac{1}{2}(1+\sigma\tau )\, =\, e_1+e_{\varphi ^*}$.}\bigskip

Les extensions biquadratiques ci-dessus sont un cas particulier 
d'{\it extensions à conjugaison complexe}. Nous entendons par là les 
corps de nombres totalement imaginaires $L$ qui sont extensions 
quadratiques d'un sur-corps totalement réel du corps de base $F$. 
Dans le contexte galoisien que nous considérons, la conjugaison 
complexe $\tau$ engendre alors un sous-groupe $\Delta_\infty$ 
d'indice 2 dans $\Delta$ (ce qui suppose donc que $\ell$ soit {\it impair}) 
et les idempotents associés $e_+=\frac{1}{2}(1+\tau )$ et $e_-=\frac{1}
{2}(1-\tau )$ permettent d'écrire tout $\Zl[\Delta]$-module $M$ comme 
somme directe de sa compo\-sante {\it réelle} $M^+ = M^{e_+}$ et de sa 
composante {\it imaginaire} $M^- = M^{e_-}$.

\bigskip

\noindent{\bf 1.b. Rappels sur le groupe des classes logarithmiques}

\medskip

Nous rappelons brièvement les principaux résultats sur les groupes 
de classes logarithmiques que nous utilisons dans la suite (cf. \cite{J6}
pour les preuves).\smallskip

Pour chaque corps de nombres $F$, notons $\J _F$ le
{\it $\ell$-adifié  du groupe  des idèles de $F$}, i.e. le produit restreint

\centerline{$\J_F\, =\, \prod_{\p }^{res}\R _{F_{\p}}$}\smallskip

\noindent des compactifiés $\ell$-adiques $\R _{F_{\p}}\, =\,
\varprojlim F_{\p} ^\times /F_{\p}^{\times \ell^n}$
des groupes multiplicatifs des complétés de $F$ aux places non complexes.

Pour chaque place finie $\p$ le sous-groupe
$\widetilde{\mathcal  U }_{F_{\p}}$
de $\R _{F_{\p} }$ formé des normes cyclotomiques (i.e. des éléments 
de $\R _{F_{\p}}$  qui sont normes à chaque étage fini de la
$\Z _{\ell}$-extension cyclotomique locale $F_{\p} ^c/F_{\p}$)
est le groupe des {\it unités logarithmiques locales} de $F_{\p}$ et le 
produit

\centerline{$\wi\U_F \, =\, \prod_{\p} \, \wi\U _{F_\p}$}\smallskip

\noindent est le groupe des {\it unités logarithmiques idéliques}~; 
c'est aussi le noyau des valuations logarithmiques
$$
\wi v_\p\ |  \ x\, \mapsto\, -\frac{\Log\, |   {x} |  _{\p}}{\deg_F\p} \ = \ 
-\frac{\Log_{Iw}\, N_{F_\p / \Q_p}   (x)}{\deg_F\p}
$$
respectivement définies sur les   $\R _{F_{\p} }$ et à valeurs
dans $\Zl$ obtenues en prenant le logarithme d' Iwasawa de la 
norme $x$ dans l'extension locale $F_{\p}/\Q _p$ corrigé par un facteur
$\deg_F\p$ choisi de sorte à assurer la surjectivité de 
$\wi v_\p$.

Le quotient 

\centerline{$\Dl _F\, =\, \J _F / \wi \U_F$}\smallskip

\noindent est, par définition, le $\ell$-groupe des
{\it diviseurs logarithmiques} de $F$~; il s'identifie, via les valuations
logarithmiques $\wi v_\p$, au $\Zl$-module libre\smallskip

\centerline{$\Dl _F\, =\, \oplus _\p \  \Zl \  \p$}\smallskip

\noindent construit sur les places finies de $F$ et on note $\wi \div_F$
la surjection naturelle de $\J_F$ sur $\Dl_F$. On définit alors le 
{\it degré} d'un diviseur logarithmique $\d  =\sum _\p n_\p \, \p$
par la formule $\deg _F \d =\, \sum _\p \ n_\p \deg_F\, \p$
et l'on écrit $\wi \Dl _F\, =\, \{ \d \in \Dl _F\, | \, \deg _F \d  \, =\, 0\, \}$ 
le sous-groupe des {\it diviseurs logarithmiques  de   degré  nul}, qui 
contient en particulier l'image $\wi \Pl _F$ dans $\Dl _F$ du sous-groupe 
de $\J _F$ formé des {\it idèles principaux}\smallskip

\centerline{$\R _F\, =\, \Zl \otimes _{\Z} F^\times$.}\smallskip

\noindent On dit que $\wi \Pl _F$ est le
sous-groupe des {\it diviseurs logarithmiques principaux}~; et
le quotient (conjecturalement fini )
$$
\wi \Cl _F\, =\, \wi \Dl _F / \wi \Pl  _F
$$
est, par définition, le {\it $\ell$-groupe des classes logarithmiques} du 
corps $F$. Enfin le noyau $\wi \E _F\, =\, \R _F \cap \ \wi \U_F$ du morphisme 
$\wi \div_F$ de $\R _F$ dans $\wi \Dl _F$ est le groupe des {\it unités 
logarithmiques globales} du corps $F$ qui apparait dans la suite exacte~:
\begin{displaymath}
\xymatrix{1 \ar[r] & \wi \E _F \ar[r] &  \R_F 
\ar[r]^{\wi \div_F} & \wi \Dl_F \ar[r] & \wi \Cl_F
\ar[r] & 1}
\end{displaymath}

Regardons maintenant les morphismes de transition qui conduisent à 
la notion, essentielle pour ce qui suit, de ramification logarithmique~:
 \smallskip

Soit donc $N/F$ une extension arbitaire de corps de nombres. Pour 
chaque nombre premier $p$, notons $\widehat{\Q}^c_p$ la 
$\widehat \Z$-extension cyclotomique de $\Q_p$, i.e. le compositum 
des $\Z_q$-extensions cyclotomiques de  $\Q_p$ pour tous les premiers 
$q$. Soit alors $\p$ un premier de $F$ au-dessus de $p$ puis ${\P}$ un 
premier de  $N$ au-dessus de $\p$. Les indices  {\it de ramification 
logarithmique} $\wi e_{N_\P /F_\p}$ et {\it d'inertie logarithmique} 
$\wi f_{N_\P/F_\p}$ sont par définition les degrés relatifs\smallskip

\centerline{$\wi e_{N_\P /F_\p} = [N_\P : N_\P \cap \widehat \Q^c_p F_\p] 
\quad  \& \quad  \wi f_{N_\P/F_\p} = [ N_{\goth P}
\cap \widehat{\mathbb Q}^c_p F_{\goth p}:F_{\goth p}]$.}\medskip

\noindent Il suit de là que l'extension, $N/F$ est logarithmiquement non
ramifiée en la place  $\P$, c'est à dire que l'on a $\wi e(N_\P /F_\p)=1$, 
si et seulement si le corps local $N_\P$ est contenu dans la $\widehat\Z$-extension
cyclotomique de $F_\p$ (i.e. si l'extension $N/F$ est localement cyclotomique 
en  $\P$). D'autre part, pour chaque $q\ne p$  les indices précédents 
pris au sens classique ou au sens logarithmique ont la même $q$-partie. 
Ils coïncident donc pour presque tout $p$ (de fait pour $p\nmid [F_\p : \Q_p]$). 

Les analogues logarithmiques des morphismes de transition sont alors 
caractérisés comme suit~: le morphisme d'extension $\wi \iota_{N/F}$ 
est donné par la formule~:\smallskip

\centerline{$\wi {\iota}_{N/F}({\goth p})=\sum_{{\goth P|\goth p}}\widetilde
 e_{N_{\goth P}/F_{\goth p}}{\goth P}$ ,}\smallskip

\noindent tandis que la norme logarithmique $\wi {N}_{N/F}$ est définie 
par l'identité ~:\smallskip
  
\centerline{$\wi N_{N/F}(\P) = \wi f_{N_\P /F_\p} \ \p$ .}\smallskip

\noindent Ces applications sont compatibles avec leurs analogues 
usuels entre  $\R_N$ and  $\R_F$ et prennent place de ce fait dans les
diagrammes commutatifs~:
$$
\CD 
\mathcal R_N     @>\wi\div_N>> \wi{\Dl}_{N}    @>\deg_N>>      \Zl @. \qquad
\qquad  @. \R_N     @>\wi\div_N>>  \wi{\Dl}_{N}    @>\deg_N>>      \Zl\\ 
@VV{N}_{N/F}V        @VV\wi { N}_{N/F}V          @\vert\text{\small \&}@.  
@AA\wi\iota_{N/F}A          @AA\wi\iota_{N/F}A       
@AA[N:F]A \\    
\R_F     @>\wi\div_{F}>> \wi\Dl_{F}   @>\deg_F>>     
\Zl @. \qquad \qquad @. \R_F     @>\wi\div_{F}>>    
\wi\Dl_{F}    @>\deg_F>>     \Zl.\\
\endCD
$$
Par passage au quotient, on obtient ainsi une description des morphismes 
de transition entre groupes de classes logarithmiques, ce qui permet 
en particulier, lorsque l'extension $N/F$ considérée est galoisienne,
d'estimer le nombre de points fixes (i.e. de classes ambiges) ou de 
copoints fixes (i.e. de classes centrales), soit~:
$$
|^G \wi \Cl_N|=|\wi \Cl_F|\ \frac{\prod\limits_{\p \in Pl^\infty_F} 
d^{ab}_\p(N/F) \ \prod\limits_{\p \in Pl^0_F} \wi e^{ab}_\p (N/F)}{[N^c : F^c] 
(\wi \E_F : \wi \E_F \cap N_{N/F}\,\R_N)}\ \kappa_{N/F}\ 
(\Dl^{I_G}_N\,\Pl_N : \wi \Dl^{I_G}_N\Pl_N),
$$
où $G$ est le groupe de Galois de l'extension~; $d^{ab}_\p(N/F)$ et 
$\wi e^{ab}_\p(N/F)$ sont res\-pectivement le degré et l'indice de 
ramification logarithmique de l'extension {\it abélienne} locale $N^{ab}
_\p\!/F_\p$ attachée à la place $\p$~; et $\kappa_{L/K}=(\mathcal N^{loc}
_{L/K} : N_{L/K}\R_L)$ désigne le nombre de noeuds de $L/K$ (cf. 
\cite {J6}), Th. 5.4). Enfin, dans le contexte galoisien défini plus haut, cette 
formule admet évidemment un raffinement composante par composante, 
qui nous sera utile plus loin (cf. \cite{JS2}, Th. 10).

\bigskip

\noindent {\bf 1.c. Introduction des noyaux étales sauvages}

\medskip

Classiquement les $\ell$-noyaux étales sauvages\footnote{ Cette appellation est 
due à T. Nguyen Quang Do \cite{N0}.} sont définis pour $i \ge 1$ 
comme les noyaux des morphismes de localisation~:\smallskip

\centerline{$\W_{2i}(F)=\Ker(H_{\text{ét}}^2(O^S_{F},\Zl(i+1))\rightarrow
 \oplus _{\p \in S}H_{\text{\'et}}^2(F_\p,\Zl(i+1)))$ ; }\smallskip

\noindent ce qui, en termes de cohomologie galoisienne, s'écrit encore~:\smallskip

\centerline{$\W_{2i}(F)=\Ker(H^2(G^S_{F},\Zl(i+1))\rightarrow 
\oplus_{\p \in S}H_{}^2({F_\p},\Zl(i+1)))$,}\smallskip

\noindent si $S$ désigne l'ensemble des places réelles ou $\ell$-adiques 
de $F$ et $\Zl(i+1)$ le $(i+1)$-ième tordu à la Tate du groupe $\Zl$. Par 
dualité de Poitou-Tate et montée dans la $\ell$-tour cyclotomique, P. 
Schneider (cf. \cite{Sc}) en tire une description des $\W_{2i}(F)$ 
comme groupe des copoints fixes d'un certain module d'Iwasawa 
tordu. En fait, cette description s'étend sans peine aux indices négatifs, 
ce qui permet de redéfinir comme suit les $\ell$-noyaux 
étales sauvages~:

\begin{Def} Soit $L$ un corps de nombres contenant les racines $2\ell
$-ièmes de l'unité, $L^c = \cup_{n \in \N} \ L[\mu_{\ell^n}]$ sa $\Zl$-extension 
cyclotomique et $\Gamma$ le groupe procyclique $\Gal(L^c/L)$. Notons 
$\Tl=\varprojlim \mu_{\ell^n}$ le module de Tate construit sur les ra\-cines 
$\ell$-primaires de l'unité, $\lT = \varprojlim \mu_{\ell^n}^*$ le module 
contagrédient, i.e. le dual de Pontrjagin de la réunion $\mu_{\ell^\infty}$
des groupes $\mu_{\ell^n}$. On définit alors le $i$-ième $\ell$-noyau 
sauvage\footnote{T. Nguyen Quang Do \cite{N2} propose la notation  
$H_{2i}L$.}de $L$ comme le quotient des copoints fixes~:\smallskip

\centerline{$\W_{2i}(L) \,\simeq\, {}^{\Gamma}(\Tl^{\otimes i}\otimes_\Zl
\C_{L^c})$,}\smallskip

\noindent où $\Tl^{\otimes i}$ est la $|i|$-ième puissance tensorielle de 
$\Tl$ pour $i \ge 0$ et de $\lT$ pour $i \le 0$, et $\C_{L^c}$ désigne le groupe de Galois $\Gal(H^{cd}_{L^c}/L^c)$ attaché à la pro-$\ell$-extension abélienne localement triviale maximale $H^{cd}_{L^c}$ de $L^c$.
\smallskip

Enfin, dans le contexte galoisien exposé au début de l'article, on 
définit le $i$-ième noyau sauvage $\W_{2i}(F)$ d'un corps de 
nombres arbitraire $F$ comme la 1-composante $\W_{2i}^{e_1}(L)$du 
noyau correspondant attaché au corps $L$.
\end{Def}

\Remarque La montée dans la tour cyclotomique $L^c/L$ ayant épuisé toute possibilité d'inertie aux places étrangères à $\ell$, le corps $H^{cd}_{L^c}$ n'est autre que la pro-$\ell$-extension abélienne maximale de $L^c$ qui est non ramifiée partout et complètement décomposée aux places au-dessus de $\ell$.

En particulier, le groupe $\C_{L^c}$ s'identifie à la limite projective \smallskip

\centerline{$\C_{L^c} \,\simeq\, \varprojlim \Cl'_{L_n}$}\smallskip

\noindent des $\ell$-groupes de $\ell$-classes $\Cl'_{L_n}$ attachés aux étages finis $L_n$ de la tour cyclotomique $L^c/L$, i.e. aux quotients des $\ell$-groupes de classes (au sens ordinaire) $\Cl_{L_n}$ par leurs sous-groupes respectifs engendrés par les classes de idéaux au-dessus de $\ell$.
\medskip

Cela étant, les groupes $\W_{2i}(F)$ admettent l'interprétation suivante~:

\begin{Prop} Soient $F$ un corps de nombres et $\ell$ un premier arbitraires.
\begin{itemize}
\item[(i)]   Pour $i > 0$, $\W_{2i}(F)$ est fini comme $i$-ème noyau sauvage de la $K$-théorie.
\item[(ii)]  Pour $i = 0$, $\ \W_{0}(F)$ est le $\ell$-groupe $\wi \Cl_F$ des classes logarithmiques de $F$.
\item[(iii)] Pour $i = -1$, $\W_{-2}(F)$ est le dual de Pontjagin du sous-groupe du radical kummérien $\Rad(L^{cd}/L^c)$ attaché à la pro-$\ell$-extension abélienne localement triviale maximale $L^{cd }$ de $L^c$ 
qui est engendré par des éléments de $F$.
\end{itemize}
\end{Prop}

\Preuve Ces divers résultats sont essentiellement bien connus~:
\begin{itemize}
\item[(i)] Pour $i>0$, l'interprétation des groupes $\W_{2i}(F)$ en termes de $K$-théorie supérieure fait intervenir les caractères de Chern étales définis par C. Soulé \cite{So} et W. Dwyer \& E. Friedlander \cite{DF}, puis considérés par G. Banaszak \cite{Ba1} et T. Nguyen Quang Do \cite{N1}. La finitude de $\W_{2i}(F)$ est ainsi un résultat profond qui prend appui sur les travaux de D. Quillen \cite{Qu} et les calculs d'A. Borel \cite{Bo}. Le cas particulier $i=1$ a été résolu par J. Tate \cite{Ta} et repose sur le théorème de finitude de Garland.

\item[(ii)] Pour $i=0$, le groupe $\ \W_{0}(F)$ n'est autre que le quotient des genres ${}^{\Gamma}\C_{F^c}$ du groupe $\C_{F^c}$ relatif à l'extension procyclique $F^c/F$, i.e. le groupe de Galois $\Gal(F^{lc}/F^c)$ attaché à la sous-extension maximale $F^{lc}$ de $H^{cd}_{F^c}$ complètement décomposée sur $F^c$ et abélienne sur $F$. Autrement dit, c'est exactement le $\ell$-groupe des classes logarithmiques défini dans \cite{J6}).

\item[(iii)] Considérons enfin le dual de Pontjagin du groupe $\W_{-2}(L)$. Nous avons~:

$\Hom(\W_{-2}(L), \Ql/ \Zl) = \Hom_\Gamma(\lT\otimes_\Zl\C_{L^c}, 
\Ql/ \Zl) = \Hom_\Gamma(\C_{L^c},\mu_{\ell^\infty})$.

\noindent Or, puisque $\C_{L^c}$ est le goupe de Galois $\Gal(L^{cd}/
L^c)$ attaché à la pro-$\ell$-extension abélienne complètement 
décomposée partout (i.e. localement triviale) maximale de $L^c$, son 
dual de Kummer $\Hom(\C_{L^c},\mu_{\ell^\infty})$ n'est autre que le 
radical associé~; et il vient, comme annoncé (cf. \cite {J3})~:

$\W_{-2}(L)^* = \Rad(L^{cd}/L^c)^\Gamma = \{\ell^{-k}\otimes x 
\in (\Ql/\Zl) \otimes L^\times |\ L^c[\sqrt[\ell^k]{x}\ ] \subset L^{cd}\}$~;

\noindent puis, en prenant les points fixes par $\Delta =\Gal(L/F)$~:

$\W_{-2}(F)^* = \Rad(L^{cd}/L^c)^{\Gamma\times\Delta} = \{\ell^{-k}\otimes x 
\in (\Ql/\Zl) \otimes F^\times |\ L^c[\sqrt[\ell^k]{x}\ ] \subset L^{cd}\}$~;

\noindent d'où le résultat annoncé.
\end{itemize}

\bigskip\medskip

\noindent{\large \bf 2. Enoncé de l'isomorphisme logarithmique et applications}

\bigskip

Notre point de départ est l'isomorphisme naturel 
de modules galoisiens qui relie les noyaux étales sauvages au groupe
des classes logarithmiques. Commençons donc par l'énoncer avec
précision dans le cadre galoisien qui nous intéresse ici~:

\begin{Th} Soit $\ell$ un nombre premier et $F$ un corps de nombres dont la $\ell$-tour cyclotomique $F[\zeta_{\ell^\infty}]/F$ est procyclique. Si $L$ est une extension abélienne de $F$ de groupe de Galois $\Delta $ d'ordre $d$ étranger à $\ell$, contenant le groupe $\mu_{\ell^r}$ des racines $\ell ^r$-ièmes de l'unité pour un $r\geq 1$, il existe pour tout $i \in \Z$ un isomorphisme canonique de $\Zl [\Delta ]$-modules\smallskip

\centerline{$(\star)\qquad {}^{\ell ^r}\!\W_{2i}(L)\, =\, \simeq \mu _{\ell ^r}^{\otimes i}
\otimes _{\Zl} \wi \Cl _L$}\smallskip

\noindent entre le quotient d'exposant $\ell ^r$ du noyau étale sauvage $\W_{2i}
(L)$ et le tensorisé $i$ fois par $\mu_{\ell ^r}$ du $\ell$-groupe des classes
logarithmiques du corps $L$.
\end{Th}

\Remarque Comme expliqué plus haut, l'isomorphisme annoncé ne 
fait intervenir {\it aucun} argument conjectural. En particulier, il ne 
présuppose vraie ni la conjecture de Gross généralisée (qui postule 
la finitude du groupe $\W_0(L) = \wi \Cl_L$), ni celle de Schneider (qui 
postule celle des groupes $\W_{2i}(F)$ pour $i<-1$), ni {\it a fortiori} celle 
de Leopoldt (qui n'intervient que pour $i=-1$)\footnote{Les liens entre 
ces diverses conjectures sont détaillés en appendice.}. 

En revanche, l'hypothèse faite $F[\zeta_{2^\infty}]=F^c$ est restrictive pour $\ell=2$.

\PreuveTh C'est un simple exercice de théorie d'Iwasawa. Introduisons 
la $\Zl$-extension cyclotomique $L^c$ de $L$~; notons $\Gamma =
\gamma^\Zl$ le groupe de Galois $\Gal(L^c/L)$~; désignons enfin par 
$\nabla_r$ l'idéal de l'algèbre d'Iwasawa $\Lambda = \Zl [[\gamma -1]]$
engendré par $\gamma -1$ et $\ell^r$. Cela étant~:

D'un côté, les résultats de Schneider (cf. \cite{Sc}) nous donnent 
l'isomorphisme~:\smallskip

\centerline{$\W_{2i}(L) \simeq {}^{\Gamma}(\Tl^{\otimes i}\otimes_\Zl
\C_{L^c})$,}\smallskip

\noindent où $\Tl^{\otimes i}$ désigne la $i$-ème puissance tensorielle 
du module de Tate $\Tl = \varprojlim \mu_{\ell^n}$ construit sur les 
racines d'ordre $\ell$-primaire de l'unité.

D'un autre côté, la théorie du corps de classes nous donne 
directement~:\smallskip

\centerline{$\wi \Cl_L\simeq {}^{\Gamma}\C_{L^c}$,}\smallskip

\noindent puisque $\wi \Cl_L$ n'est autre que le quotient des genres 
associé à $\C_{L^c}$ relativement à l'extension procycique $L^c/L$ 
(cf. \cite{J6}). Il vient donc~:
\begin{align*}
^{\ell^r}\!\W_{2i}(L) &\simeq (\Tl^{\otimes i}\otimes_\Zl \C_{L^c}) / \nabla_r
(\Tl^{\otimes i}\otimes_\Zl \C_{L^c}) \\
&\simeq \Tl^{\otimes i}\otimes_\Zl (\C_{L^c} / \nabla_r \C_{L^c}) 
\simeq \Tl^{\otimes i}\otimes_\Zl {}^{\ell^r}\! \Cl_L \simeq
\mu_{\ell^r}^{\otimes i} \otimes _\Zl  \wi \Cl_{L},
\end{align*}
comme annoncé, puisque par hypothèse le groupe $\Gamma$ opère 
trivialement sur $\mu_{\ell^r}$.

\bigskip

\noindent{\bf 2.a. Formules de rang et surjectivité de la descente}

\medskip

Commençons par énoncer le Théorème 3 pour chaque composante isotypique.

\begin{Th} Soit $F$ arbitaire\footnote{Sous la seule condition de 
procyclicité de sa 2-extension cyclotomique dans le cas $\ell=2$.} et $L$ 
comme plus haut. Pour chaque idempotent  $e_{\varphi}$ de l'algèbre 
 $\Zl[\Delta ]$ et  $i \in \Z$, il existe un isomorphisme de $\Z_\phi$-modules~:

\centerline{${}^{\ell ^r}\!\W_{2i}^{e_{\varphi \bar\omega^i}}\!(L)\, \simeq\,
  {}^{\ell ^r} \wi \Cl^{e_\varphi}_L$,}\smallskip

\noindent entre la $\varphi\bar\omega^i$-composante du quotient 
d'exposant $\ell^r$ du noyau sauvage $\W_{2i}(L)$ et la  
$\varphi$-composante du quotient d'exposant $\ell^r$ du groupe des 
classes logarithmiques.\smallskip

En particulier, on a~: \quad $^{\ell ^r}\!\W_{2i}(F)\, \simeq\,
  {}^{\ell ^r} \wi \Cl^{e_{\varphi\omega^i}}_L$.
\end{Th}

\Preuve C'est la traduction directe du Théorème 3 composante par 
composante\footnote{On retrouve ainsi en particulier les résultats de 
périodicité obtenus par M. Kolster \cite{K3}.}.\medskip

Dans le cas particulier des corps biquadratiques $L = \Q [\sqrt{-3},
\sqrt d]$ étudié par Scholz, le Théorème ci-dessus prend une forme particulièrement simple~:

\begin{Cor} Soit $k=\Q[\sqrt d]$ un corps quadratique (réel ou 
imaginaire, mais distinct de $\Q[\sqrt{-3}]$)
et $k^*=\Q[\sqrt{-3d}]$ son reflet dans l'involution du miroir. Alors pour
chaque ${i>0}$, les quotients d'exposant 3 des noyaux sauvages étales 
de $k$ sont donnés à partir des 3-goupes de classes logarithmiques 
de $k$ et $k^*$ par~:
$$
{}^3\W_{2i}(k) \simeq \left\{ \begin{array}{ll}
\mu_3^{\otimes i} \otimes_{\Z_3} \wi \Cl_k \quad & pour \  i \ pair, \\
\mu_3^{\otimes i} \otimes_{\Z_3} \wi \Cl_{k^*} \quad & pour \ i \ impair.
\end{array} \right.
$$
\end{Cor}

\Remarque Dans le cas  du corps  quadratique
$k = \Q[\sqrt{-3}]$, il vient directement~:\smallskip

\centerline{${}^3\W_{2i}(k) \simeq \mu_3^{\otimes i} \otimes_{\Z_3} 
\wi \Cl_k = 1$,}\smallskip

\noindent puisque le corps cyclotomique $\Q[j]$ est 3-régulier 
donc 3-logarithmiquement principal (cf. \cite{GJ}). Tous les $\W_{2i}(k)$ 
ont donc dans ce cas une 3-partie triviale.\medskip

A l'instar de ce qui est fait dans \cite{JS2}, l'isomorphisme $(\star)$ 
permet de transporter aux 
noyaux étales sauvages les inégalités du miroir sur les $\ell$-rangs des 
classes logarithmiques, qui ne sont en fin de compte que la transcription
dans le cadre logarithmique du très classique {\it Spiegelungssatz} de
Leopoldt~: 

\begin{Prop}Supposons $F$ totalement réel et $L$ à conjugaison 
complexe i.e. extension quadratique totalement imaginaire d'un 
sur-corps totalement réel de $F$. Alors, pour chaque caractère $\ell
$-adique irréductible imaginaire $\varphi $ du groupe $\Delta $, on a les 
inégalités entre $F_\phi$-dimensions des composantes isotypiques~:
\smallskip

\centerline{$ 0 \leq \dim_{F_\phi} {}^\ell \W_{2i}^{e_{\varphi\bar\omega^i}}
\!(L) - \dim_{F_\phi} {}^\ell \W_{2i}^{e_{\varphi^*\bar\omega^i}}\!(L) \leq \ [F:\Q ]$.}
\end{Prop}

\Preuve D'après l'isomorphisme ($\star$), la différence considérée 
s'écrit encore~:\smallskip

\centerline{$\dim_{F_\phi} {}^\ell \wi \Cl _L^{e_\varphi}  -
\dim_{F_\phi} {}^\ell \wi \Cl_L^{e_{\varphi^*}}$.}\smallskip

\noindent Le résultat annoncé résulte donc directement du Théorème 
4 de \cite{JS2}.

\begin{Cor} Prenons $F=\Q$ et $L=\Q[\zeta _\ell]$. Le groupe de 
Galois $\Delta$ est alors cyclique d'ordre $\ell-1$ et ses caractères 
$\ell$-adiques irréductibles $\varphi$ sont de dimension 1. Pour 
chaque $\varphi$ réel, on a donc l'inégalité entre $\ell$-rangs des 
noyaux sauvages~:\smallskip

\centerline{$ 0 \ \leq \ \rg_\ell  {}^\ell \W_{2i}^{e_{\varphi ^*\bar\omega^i}}
\!(L)\ - \ \rg_\ell {}^\ell \W_{2i}^{e_{\varphi\bar\omega^i }}\!(L)\ \leq 1$.}
\end{Cor}

\begin{Cor} Prenons $F=\Q$ et $L=\Q [\sqrt d,\sqrt {-3}]$. Notons 
$k=\Q \ [\sqrt d]$ le sous-corps quadratique réel de $L$ et $k^* = 
\Q[\sqrt {-3d} ]$ son reflet. Il vient alors~:\smallskip

\centerline{$ 0 \leq  \rg_3 {}^3 \W_{4i}(k)  -  \rg_3 {}^3 \W_{4i}(k^*)  
=  \rg_3 {}^3 \W_{4i+2}(k^*)  -  \rg_3 {}^3 \W_{4i+2}(k) \leq 1$.}
\end{Cor}

Notons d'autre part qu'il est possible, sous certaines conditions de 
déduire de l'isomorphisme ($\star$) des résultats plus forts sur les 
groupes tout entiers~:

\begin{Prop} Conservons les notations du Théorème 4 et supposons 
que pour un caractère $\ell$-adique irréductible $\varphi$ du groupe 
$\Delta$, la $\varphi$-composante du $\ell$-groupe des classes 
logarithmiques $\wi \Cl_L$ soit finie, disons d'ordre $\ell^{s_\varphi}$. 
Alors~:\begin{itemize}

\item[(i)] Si l'on a $s_\varphi<r$, la $\varphi\bar\omega^i$-composante 
de $\W_{2i}(L)$ est encore finie et on a~:\smallskip

\centerline{$\W_{2i}^{e_{\varphi\bar\omega^i}}\!(L) \ \simeq \ 
\wi \Cl^{e_\varphi}_L$.}

\item[(ii)] Si l'on a $s_\varphi \le r$, la même conclusion vaut encore sous
réserve d'injectivité~\footnote{Cette injectivité sera discutée plus loin 
en liaison avec l'étude des noyaux de capitulation.} 
du morphisme d'extension~: $\W_{2i}^{e_{\varphi\bar\omega^i}}\!(L) 
\rightarrow \W_{2i}^{e_{\varphi\bar\omega^i}}\!(L[\zeta_{\ell^{r+1}}])$.
\end{itemize}
\end{Prop}

\Preuve Identique {\it mutatis mutandis} à celle du Théorème 7 de 
\cite{JS2}.\medskip

Expliquons enfin comment l'isomorphisme ($\star$) permet de résoudre
très simplement la question de la surjectivité du morphisme de descente
(ou corestriction) dans une extension arbitraire de corps de nombres (cf. \cite{Ka})~:

\begin{Prop} Soit $F$ un corps de nombres arbitraire et $N$ une 
extension de $F$ ayant pour degré une puissance de $\ell$ (mais non 
nécessairement galoisienne). Soit $L$ abélienne sur $F$, de groupe 
$\Delta$, de degré relatif $d$ étranger avec $\ell$, contenant les 
racines $2\ell$-ièmes de l'unité, et $LN$ le compositum de $L$ et de 
$N$. Alors, pour tout $i \in \Z$ l'application
canonique de corestriction $Tr_{N/F}$\smallskip

\centerline{$\W_{2i}(N) \rightarrow \W_{2i}(F)$}\smallskip

\noindent est surjective si et seulement si LN/L ne contient pas de 
sous-extension cyclique de degré $\ell$ et isotypique\footnote{ C'est à 
dire dont le groupe de Galois est un $\Zl[\Delta]$-module isotypique.} 
de caractère $\omega^i$ qui soit logarithmiquement non ramifiée et 
disjointe de $L^c/L$. En particulier, lorsque l'extension $N/F$ considérée 
est galoisienne, on a toujours $Tr_{N/F}(\W_{2i}(N)) = \W_{2i}(F)$ dès qu'on a $\,\omega^i \ne 1$.
\end{Prop}

\Preuve  Prenant, en effet, les 1-composantes, nous obtenons 
immédiatement~:\smallskip

\centerline{$\W_{2i}(F)/Tr_{N/F}(\W_{2i}(N))\W_{2i}(F)^\ell \simeq 
\wi\Cl^{e_{\omega^i}}_L / \wi N_{LN/L}(\wi\Cl^{e_{\omega^i}}_{LN})
\wi\Cl^{e_{\omega^i}\, \ell}_L$.}\smallskip

\noindent Et la théorie du corps de classes (cf. \cite{J7}) 
nous donne bien le résultat annoncé, en vertu de l'interprétation de 
$\wi \Cl_{L}$ comme groupe de Galois $\Gal(L^{lc}/L^c)$ attaché à la 
$\ell$-extension abélienne logarithmiquement non ramifiée maximale 
de $L$.

\bigskip

\noindent{\bf 2.b. Trivialité des noyaux étales supérieurs}

\medskip

Une autre application essentielle de l'isomorphisme $(\star)$ est 
de caractériser en termes de classes logarithmiques la trivialité des 
diverses composantes iso\-typiques des $\ell$-noyaux étales sauvages, 
puisque, dans le contexte galoisien précisé plus haut, on a évidemment 
l'équivalence~:\smallskip

\centerline{$\W_{2i}^{e_{\varphi \bar\omega^i}}\!(L) = 1\ \Leftrightarrow \ 
{}^\ell \W_{2i}^{e_{\varphi \bar\omega^i}}\!(L) = 1\ \Leftrightarrow \
{}^\ell \wi \Cl^{e_\varphi}_L = 1\ \Leftrightarrow \ \wi \Cl^{e_\varphi}_L = 1$.}
\smallskip

\noindent En particulier, on dispose là d'un moyen d'étudier la 
trivialité de $\W_{2i}(F)$ pour n'importe\footnote{Toujours sous la seule 
condition de procyclicité de la 2-extension cyclotomique pour $\ell=2$.}
quel corps de nombres $F$, et notamment la propagation de cette 
trivialité dans une $\ell$-extension de tels corps, en transportant par 
l'isomorphisme $(\star)$ la formule des classes logarithmiques 
centrales (cf. \cite{JS2}). Il vient ainsi~:

\begin{Th} Soit $N/F$ une $\ell$-extension de corps de nombres de 
groupe de Galois $G$ et $LN/F$ l'extension obtenue par composition 
avec une extension abélienne $L/F$, de groupe $\Delta$, de degré $d$
relatif étranger avec $\ell$, contenant les racines $2\ell$-ièmes de 
l'unité. Alors, pour tout caractère $\ell$-adique irréductible $\varphi 
\ne 1$ du groupe $\Delta$, on a l'équivalence~:
\begin{equation*}
\W_{2i}^{e_{\varphi \bar\omega^i}}\!(LN) = 1 \Leftrightarrow \begin{cases}
\W_{2i}^{e_{\varphi \bar\omega^i}}\!(L) =  1 \qquad { \rm \&  }\\
(\wi \E_{LN}^{e_\varphi} :
\wi \E_{LN}^{e_\varphi} \cap \mathcal  N_{LN/L}^{loc})  = 
\prod_{\p }\, \wi e_\p^{ab}(N/F) ^{<\varphi ,\chi _{\p} >},
\end{cases}
\end{equation*}
où $\wi e_\p^{ab}(LN/L)$ désigne l'indice de ramification logarithmique 
abélianisé de laplace $\p$ de $F$ dans l'extension $N/F$et $\chi_\p$ 
l'induit à $\Delta$ du caractère unité de son sous-groupe de 
décomposition $\Delta_\p$ dans $L/F$, tandis que $\wi \E_{LN}$ est le 
groupe des unités logarithmiques du corps $LN$ et $\mathcal  N_{LN/L}
^{loc}$ le groupe des normes locales dans $NL/L$.
\end{Th}

\Preuve Il suffit de transporter par l'isomorphisme $(\star)$ le Théorème
10 de \cite{JS2}.\medskip

\Remarque Outre des simplifications techniques permettant de tuer la 
contribution éventuelle du groupes des noeuds relatif à l'extension  
$LN/L$, l'hypothèse $\varphi \ne 1$ assure essentiellement la validité 
de la descente~:\smallskip

\centerline{$\W_{2i}^{e_{\varphi \bar\omega^i}}\!(LN) = 1 \ \Rightarrow \ 
\W_{2i}^{e_{\varphi \bar\omega^i}}\!(L) = 1$.}\smallskip

\noindent Sans surprise, la montée tout au contraire met en jeu de 
façon non triviale l'arithmétique (logarithmique) de l'extension 
considérée.\par

On observera que pour $\varphi = 1$, en revanche,  la descente  peut 
être en défaut~: c'est en particulier le cas lorsque le corps de base $F$ 
n'est pas logarithmiquement principal mais possède une $\ell$-tour 
localement cyclotomique finie, i.e. une $\ell$-extension (finie) qui est 
logarithmiquement principale (cf. \cite{JS1}).

\begin{Sco} Les conclusions du Théorème valent encore pour $\varphi
 = 1$, dès lors qu'une au moins des places de $F$ se ramifie totalement 
(au sens logarithmique) dans l'extension $N/F$.
\end{Sco}

\Preuve Cette hypothèse anéantit, en effet, tous les facteurs parasites
(cf. \cite{J6}.

\begin{Cor} Supposons que l'extension $L/F$ soit à conjugaison 
complexe. Alors,  pour tout caractère imaginaire irréductible 
$\varphi \ne \omega$, on a l'équivalence~:\smallskip

\centerline{$\W_{2i}^{e_{\varphi \bar\omega^i}}\!(LN) = 1 \ \Leftrightarrow \ 
\W_{2i}^{e_{\varphi \bar\omega^i}}\!(L) = 1 \ \& \ <\varphi ,\chi _{\p} >=0 
\quad \forall \ \p\in \wi R_{N/F}$,}\smallskip

\noindent où $\wi R_{N/F}$ est l'ensemble des places de $F$ logarithmiquement
ramifiées dans $N/F$.  

Et la même équivalence vaut pour $\varphi\ = \omega$, à ceci près que 
l'on peut avoir dans ce cas $<\varphi ,\chi _{\p} > \ne 0$ en une place 
$\p_\circ$ de $\wi R_{N/K}$ au plus, qui vérifie en outre~:\smallskip

\centerline{$(\mu_L : \mu_L \cap \mathcal  N^{\rm loc}_{LN/L}) \ = \ 
\wi e_{\p _\circ }^{ab}(N/F)^{<\omega ,\chi_\p>}$.}
\end{Cor}

\Preuve
Il suffit d'observer que les seules unités logarithmiques  imaginaires
sont les racines de l'unité, qui forment un $\Zl [\Delta ]$-module
isotypique de caractère $\omega$.

\begin{Cor} Supposons $\ell$ impair, $F$ totalement réel et prenons 
$L=F[\zeta_\ell]$. Soit $N$ une $\ell$-extension (galoisienne) 
arbitaire de $F$. Alors~:

\item[(i)] Pour $i \not\equiv 1\ [{\rm mod }\ d]$ impair, on a $\W_{2i}
(N) =1$ si et seulement si les deux conditions suivantes sont réunies~:

(i,a) $ \W_{2i}(F) = 1$ et
 
(i,b) les  places $\p$ de $F$ logarithmiquement ramifiées dans $N/F$
ne sont pas complètement décomposées dans la sous-extension 
$L^{\Ker \omega^i}/F$ .

\item[(ii)] Pour $i \equiv 1\ [{\rm mod }\ d]$, ces deux conditions
doivent être remplacées par~:

(ii,a) $ \W_{2i}(F) = 1$ et

(ii,b) Une place au plus $\p_\circ$ de $F$ est à la fois logarithmiquement 
ramifiée dans $N/F$ et complètement décomposée dans $L/F$~; de plus on 
a alors~:\smallskip

\centerline{$\wi e_{\p_\circ} ^{ab}(N/F)\, =
\, (\mu_L:\mu_L \cap \mathcal  N_{LN/L}^{loc})$.}

\end{Cor}

\Preuve C'est immédiat d'après le corollaire 13 appliqué avec $\varphi 
=\omega^i$, la condition d'orthogonalité $<\omega^i , \chi _\p ,   >\, =\, 0$
exprimant simplement la non trivialité de $\omega^i$ sur le sous-groupe  
de décomposition $\Delta _\p$ de la place $\p$ dans l'extension $L/F$.\medskip

A l'opposé, la Théorie des Genres  fournit également une minoration 
du rang~:

\begin{Prop} Soit $F$ un corps de nombres arbitraire et $\ell$ un nombre
premier. Alors $F$ possède une infinité de $\ell$-extensions cycliques 
$N$ de degré $\ell$ dont tous les noyaux sauvages $\W_{2i}(N)$ ont un 
$\ell$-rang arbitrairement grand.
\end{Prop} 

\Preuve Il suffit de transporter par ($\star$) les
 minorations du rang logarithmique données par la Théorie des 
Genres, dont le principe est le suivant~: La suite exacte des classes 
logarithmiques centrales (cf. \cite{J6}) appliquée dans une $\ell$-extension 
cyclique $N/F$ de corps de nombres fait apparaître plusieurs termes 
dont {\it tous} ont une contribution explicitement bornée, à l'exclusion du 
terme de ramification d'ordre $\prod_\p \wi e_{N_\p/F_\p}$ qui est arbitrairement 
grand avec le nombre de ramifiés. Concrètement, un calcul sans surprise 
mené composante par composante dans l'extension composée $NL/L$
avec $L=F[\zeta_\ell]$ donne\footnote{L'inégalité vaut pour $\ell$ 
impair ; et à 1 près pour $\ell = 2$, ce qui est donc sans conséquence.}~:
\smallskip

\centerline{${\rm rg}_\ell \W_{2i}(N) \, = \, {\rm rg}_\ell \wi \Cl_{NL}^{e_{\omega^i}} \, 
\ge \, {\rm rg}_\ell {}^G \wi \Cl_{NL}^{e_{\omega^i}}\,  \ge  \, <\omega^i , \sum_
{\p \in \wi R_{N/F}} \chi_\p - \sum_{\p | \infty}\chi_\p - 1>$,}\smallskip

\noindent d'où le résultat, le caractère $\sum_{\p \in \wi R_{N/F}} 
\chi_\p$ pouvant être pris arbitrairement grand avec le nombre de places 
ramifiées (au sens ordinaire ou logarithmique).

\bigskip

\noindent{\bf 2.c. Le problème de la capitulation}

\medskip 

Venons en maintenant au noyau du morphisme d'extension déjà étudié par 
J. Coates, R. Greenberg, B. Kahn, T. Nguyen quang Do ou l'auteur \cite{J4} dans divers contextes. Du point de vue de la Théorie d'Iwasawa, la question de la 
capitulation pour les noyaux étales sauvages dans une $\Zl$-extension se 
présente comme suit~:

On dispose d'un module noethérien $X$ sur l'algèbre d'Iwasawa $\Lambda 
= \Zl [[\gamma -1]]$ attachée à un groupe procyclique $\Gamma = \gamma
^\Zl$~; on note $\omega_n =\gamma^{\ell^n} -1$~; et on s'intéresse aux 
noyaux des morphismes de transition $Cap_n = \Ker (X_n \mapsto X_m)$ 
pour $m \gg n \gg 0$ induits par la multiplication par $\omega_m / 
\omega_n$ entre les quotients $X_n = X/\omega_n X$ et $X_m = X/
\omega_m X$ pour $n$ et $m-n$ assez grands~; ce qui revient à 
considérer le noyau $Cap_n = \Ker (X_n \mapsto X_\infty)$ du 
morphisme d'extension à valeurs dans la limite inductive $X_\infty$ des 
$X_n$. Notant $T$ le plus grand $\Lambda$-sous-module fini\footnote{Autrement dit, 
le noyau des morphismes de localisation $X \mapsto X_\wp$, où $\wp$ décrit 
l'ensemble des idéaux premiers de hauteur 1 de l'anneau $\Lambda$.} de 
$X$, on a, par un calcul élémentaire~:\smallskip

\centerline{$Cap_n = \{ x+\omega_n X \ | \ (\omega_m / \omega_n)\ x  =
0 \} \underset{m \gg n}{=} (T + \omega_n X ) / \omega_n X \underset{n \gg 0}
{\simeq} T$,}\smallskip

\noindent où l'on voit que non seulement les groupes $Cap_n$ sont 
bornés (comme l'a montré Iwasawa pour les groupes de classes) mais 
qu'ils sont ultimement isomorphes et que leur limite projective a une 
interprétation très simple. Lorsque, de plus, le groupe $X$ est un $\Zl
$-module noethérien\footnote{C'est à dire un $\Lambda$-module de 
torsion dont l'invariant {\it mu} d'Iwasawa est trivial.}, ce qui est la 
situation standard pour la limite projective des $\ell$-groupes de classes 
d'idéaux dans la $\Zl$-extension cyclotomique d'un corps absolument 
abélien, d'après un résultat de Ferrero et Washington conjecturé 
pour tout corps de nombres (cf. \cite{Wa}), on peut même dire mieux~: d'après \cite{JG}, 
en effet, il existe alors une famille $(\alpha_j)_{j=1,\dots,\lambda}$ d'entiers 
relatifs tels qu'on ait, pour tout $n$ assez grand~:\smallskip

\centerline{$X_n \simeq \big( \oplus_{j=1}^\lambda \ \Z/\ell^{n+\alpha_j}\Z
\big) \ \oplus \ Cap_n$~;}\smallskip

\noindent de sorte que $Cap_n \simeq T$ est alors un facteur direct 
du $\Zl$-module $X_n$. \smallskip

Enfin, dans le contexte galoisien qui nous intéresse, les groupes $X$ 
que nous considérons sont naturellement des $\Lambda [\Delta]$-modules, 
ce qui permet de généraliser les résultats ci-dessus aux diverses 
composantes isotypiques regardées comme $\Z_\phi[[\gamma-1]]
$-modules, ainsi qu'il est fait dans \cite{J2}. Il vient ainsi~:

\begin{Th} Soit $F$ un corps de nombres arbitraire et $L$ une extension 
abélienne de $F$, contenant les racines $2\ell$-ièmes de l'unité, de 
groupe de Galois $\Delta$, de degré relatif $d$ étranger avec $\ell$, 
puis $F_\infty = \cup_{n\in\N} F_n$ la $\Zl$-extension cyclotomique de 
$F$ et $L_\infty = \cup_{n\in\N} L_n$ celle de $L$. Notons $\Lambda [
\Delta] = \Zl[[\gamma-1]][\Delta]$ l'algèbre d'Iwasawa construite sur un 
générateur topologique $\gamma$ du groupe procyclique $\Gamma =
 \Gal(L_\infty/L) \simeq \Gal(F_\infty /F)$ et $\T_{L_\infty}$ le $\Lambda 
[\Delta]$-sous-module fini du groupe de Galois $\ \C_{L_\infty} = 
\varprojlim \wi \Cl_{L_n}$ de la $\ell$-extension abélienne complètement 
décomposée maximale de $L_
\infty$. Avec ces notations~:\begin{itemize}

\item[(i)] Le $\Zl[\Delta]$-module $\T_{L_\infty}$ est la limite projective 
$\varprojlim \wi {Cap}_{L_n}$ des groupes de capitulation logarithmique 
$\wi {Cap}_{L_n} = \Ker\ (\wi\Cl_{L_n} \mapsto 
\wi\Cl_{L_\infty})$. 

\item[(ii)] Les noyaux de capitulation $C\W_{2i}(L_n) =
\Ker\ (\W_{2i}(L_n) \mapsto \W_{2i}(L_\infty))$ sont donnés pour tout 
$n \gg 0$ et tout $i \in \Z$ par l'isomorphisme galoisien~:\smallskip

\centerline{ $C\W_{2i}(L_n) \ \simeq \ \Tl^{\otimes i} \otimes \T_{L_\infty}
\ \simeq \ \Tl^{\otimes i} \otimes \wi {Cap}_{L_n}$.} 

\end{itemize}
\end{Th}

\begin{Sco} Si en outre le groupe $\C_{L_\infty}$ est de type fini sur $\Zl$
(i.e. si l'invariant mu d'Iwasawa attaché aux $\ell$-groupes de classes 
au sens ordinaire dans la tour $L_\infty/L$ est nul), le sous-module 
$C\W_{2i}(L_n)$ est un facteur direct du $\Zl[\Delta]$-module $\W_{2i}(L_n)$
pour tout $n$ assez grand~; et on a pour tous $r \ge 1$ et $n \gg 0$~:
\smallskip

\centerline{$^{\ell^r}\!\W_{2i}^{e_{\varphi}}\!(L_n) \ \simeq \ 
\big( \Z_\phi / \ell^r \Z_\phi \big)^{\lambda_{\varphi\bar\omega^i}} \ \oplus \ 
{}^{\ell^r}\!\T_{L_n}^{e_{\varphi\bar\omega^i}}$,}\smallskip

\noindent où $\lambda_{\varphi\bar\omega^i}$ est l'invariant lambda 
d'Iwasawa\footnote{C'est à dire la dimension sur le corps des fractions
$\Q_\phi$ de $\Z_\phi$ du  $\Q_\phi$-espace $\Q_\phi \otimes_{\Z_\phi} 
\C_{L_\infty}^{\varphi\bar\omega^i}$}attaché au $\Z_\phi[[\gamma-1]]
$-module $\C_{L_\infty}^{\varphi\bar\omega^i}$.
\end{Sco}

\Preuve L'assertion {\it (i)} s'obtient en appliquant ce qui précède au 
$\Lambda[\Delta]$-module $\C_{L_\infty}$~; l'assertion {\it (ii)} en 
l'appliquant au module tordu $\Tl^{\otimes i} \otimes \C_{L_\infty}$ dont 
le plus grand sous-module fini est évidemment $\Tl^{\otimes i} \otimes 
\T_{L_\infty}$. Enfin le Scolie s'ensuit par spécialisation aux $\varphi
\bar\omega^i$-composantes, sous réserve de trivialité de la $\varphi
\bar\omega^i$-partie $\mu_{\varphi\bar\omega^i}$ de l'invariant {\it mu} 
d'Iwasawa.

\Remarque Comme établi dans \cite{J2} (et généralisé dans \cite{JMa}),
le groupe $\C_{L_\infty}$ (et plus généralement tous les $\ell$-groupes 
$\C_T^S(L_\infty) = \varprojlim \Cl^S_T(L_n)$ de $S$-classes $T$-infinitésimales) 
ont même invariant\footnote{Il s'agit donc de l'invariant {\it mu} 
d'Iwasawa habituel.} {\it mu} d'Iwasawa que le groupe $\varprojlim \Cl_{L_n}$.

\begin{Cor} Avec les notations du Théorème, la $\varphi$-partie du 
noyau de capitulation sauvage est donné pour $n \gg 0$, à partir de la 
cohomologie des unités logarithmiques, par l'isomorphisme de 
$\Z_\phi$-modules~:\smallskip

\centerline{$C\W_{2i}^{e_\varphi}\!(L_n) \ \simeq \ H^1(\Gamma_n, \wi\E
^{e_{\varphi\omega^i}}_{L_\infty})$.}\smallskip

\noindent En particulier, lorsque l'extension $L/F$ est à conjugaison 
complexe, les groupes étales réels $C\W_{4i+2}(F)$ et les groupes 
imaginaires $C\W_{4i}^-(L)$ sont triviaux.
\end{Cor}

\Remarque On retrouve là le résultat de trivialité de Kolster et 
Movahhedi \cite{KM}.

\Preuve Les unités logarithmiques imaginaires se réduisent, en effet, 
aux racines de l'unité, lesquelles sont cohomologiquement triviales 
dans la tour $L_\infty/L$.

\begin{Cor} Lorsque l'extension $L/F$ est à conjugaison complexe, et 
sous les conjectures d'Iwasawa et de Greenberg dans la tour $L_\infty
/L$, aucun des groupes $C\W_{4i}(F_n)$ ou $C\W_{4i+2}^-(L_n)$ ne peut être 
ultimement trivial sans que tous les groupes $\W_{4i}(F_n)$ et $\W_
{4i+2}^-(L_n)$ ne le soient également.
\end{Cor}

\Preuve Rappelons (cf. appendice) que sous les conjectures d'Iwasawa 
et de Greenberg, l'invariant {\it mu} et la partie plus de l'invariant {\it 
lambda} associés au groupe $\C_{L_\infty}$ sont nuls~;  les groupes 
$\wi \Cl^{e_\varphi}_{L_n}$ attachés aux caractères réels sont alors 
d'ordre borné~; et d'après le Scolie il en va de même des groupes 
$\W_{2i}^{e_{\varphi\omega^i}}\!(L_n)$ qui coïncident donc avec leurs
sous-groupes de capitulation respectifs $C\W_{2i}^{e_{\varphi
\omega^i}}\!(L_n)$.

\newpage

\centerline{\large \bf Appendice} \medskip

\centerline{\bf Les conjectures standard de la Théorie d'Iwasawa}

\bigskip

La  principale conjecture sur les noyaux étales sauvages s'énonce 
comme suit~:

\begin{Conj} Les groupes $W_{2i}(F)$ sont finis pour tout corps de 
nombres $F$.
\end{Conj}

\Remarques De fait la conjecture ci-dessus réunit plusieurs résultats 
classiques~:

{\it (i)}   Pour $i > 1$, la finitude des noyaux sauvages $\W_{2i}(F)$, 
est un résultat profond de $K$-théorie supérieure, comme expliqué plus haut dans la section 1.c (cf. \cite{K4}).

{\it (ii)}   Pour $i =1$, la finitude du noyau sauvage $\W_{2}(F)$, 
naguère {\it conjecture principale}, est une conséquence
d'un théorème de Garland (cf. \cite{Ta}, \cite{Sc}).

{\it (iii)}  Pour $i = 0$, la finitude de $\W_{0}(F)= \wi \Cl_F$ n'est rien 
d'autre que la {\it conjecture de Gross généralisée} (cf. \cite{J6}).

{\it (iv)} Pour $i = -1$, la finitude de $\W_{-1}(F)$ est une forme
équivalente mais peu connue de la {\it conjecture de Leopoldt} pour le 
corps $F$.

{\it (v)} Pour $i < -1$ enfin, la finitude des noyaux $\W_{2i}(F)$ constitue  
à proprement parler la {\it conjecture de Schneider} (cf. \cite{Sc}).

\medskip\noindent
{\it Preuve de (iv).} D'après \cite{J7}, une façon d'énoncer la {\it conjecture de 
Leopoldt} dans le corps $F$ consiste à écrire la trivialité du {\it groupe
des unités infinitésimales}~:\smallskip

\centerline{$\E_F^\infty = \{ \varepsilon \in \Zl \otimes_\Z E_F | \ 
\varepsilon_\l = 1 \ \forall \l \mid \ell\}$,}\smallskip

\noindent qui en constitue précisément le groupe de défaut.
Cela étant, puisque le groupe $\W_{-1}(F)$ est de type fini, 
son dual de Pontjagin $\W_{-1}(F)^*$ est de cotype fini et affirmer la 
finitude de $\W_{-1}(F)$ revient à postuler que le sous-module divisible 
maximal de $\W_{-1}(F)^*$ est trivial. L'équivalence résulte donc du lemme~:

\begin{Lem} On a : $(\W_{-1}(F)^*)^{\rm div} = (\Ql/\Zl) 
\otimes_\Zl \E_F^\infty$.
\end{Lem}

\Preuve Prenons $L= F[\zeta_{2\ell}]$ comme plus haut et observons 
que si $\ell^{-k} \otimes x$ est  divisible dans $\W_{-1}(F)^*$, l'extension 
procyclique $L^c[\sqrt[\ell^\infty]{x}\ ]/L^c$ n'est localement triviale 
partout que si $x$ est une 
unité en dehors de $\ell$ et une racine de l'unité aux places divisant 
$\ell$, c'est à dire finalement que pour $\ell^{-k} \otimes x \in (\Ql/\Zl) 
\otimes_\Zl \E_F^\infty$.\medskip

De la description de Schneider $\W_{2i}(L) \simeq {}^{\Gamma}(\Tl^{\otimes i}\otimes_\Zl
\C_{L^c})$, il suit~:

\begin{Sco} Le corps $F$ étant donné, soit $L= F[\zeta_{2\ell}]$ comme 
plus haut. Fixons un générateur topologique 
$\gamma$ du groupe $\Gamma = \Gal(L^c/L) \simeq \Gal(F^c/F)$~; notons 
$\Lambda[\Delta]=\oplus\Lambda_\varphi$ l'algèbre d'Iwasawa $\Zl
[\Delta][[\gamma-1]]$ et $\kappa$ le caractère de l'action de 
$\Gamma$ sur les racines de l'unité. Avec ces notations, la finitude 
du groupe $\W^{e_\varphi}_{2i}\!(L)$ traduit directement le fait que 
$\kappa^{-i}(\gamma)$ n'est pas racine du polynôme caractéristique 
$P_{\varphi\bar\omega^i}(\gamma) \in \Z_\phi[\gamma]$ attaché au 
$\Lambda_\phi$-module $\ \C^{e_{\varphi\bar\omega^i}}_L$.
\end{Sco}

La {\it conjecture de Greenberg} sur les classes des corps totalement 
réels est, de fait,  beaucoup plus radicale pour ce qui est des composantes 
réelles~:

\begin{Conj} Soit $F$ totalement réel et $L= F[\zeta_{2\ell}]$ comme 
plus haut. On a alors~: $\deg P_{\varphi
\bar\omega^i}(\gamma) =0$, i.e. $P_{\varphi\bar\omega^i}(\gamma) \in 
\Z_\phi$, dès que le caractère $\varphi\bar\omega^i$ est réel. 
\end{Conj}

\Preuve Rappelons que la conjecture de Greenberg postule la nullité 
de l'invariant {\it lambda} d'Iwasawa attaché à la limite projective 
$\varprojlim \Cl_{F_n}$ des $\ell$-groupes de classes au sens ordinaire 
dans la $\Zl$-extension $F_\infty/F$ d'un corps de nombres totalement 
réel ou, ce qui revient au même, de la partie plus de l'invariant {\it 
lambda} du groupe $\varprojlim \Cl_{L_n}$ pour un corps $L$ à 
conjugaison complexe. {\it A fortiori} donc en 
est-il de même lorsqu'on remplace la limite projective précédente par 
son quotient $\C_{L_\infty}$. D'où la conséquence énoncée, l'invariant
{\it lambda} d'un $\Lambda$-module n'étant autre que le degré de son 
polynôme caractéristique.

\newpage

Enfin la {\it conjecture d'Iwasawa}\footnote{Pour $F$ abélien sur $\Q$, c'est le théorème de Ferrero \& Washington (cf. \cite{Wa}).} sur les invariants {\it mu} (cf. \cite{Iw}) s'écrit~:

\begin{Conj} Le groupe $\C_{F_\infty}$ est de type fini sur $\Zl$ pour tout $F$~; 
en d'autres termes, son invariant {\it mu} est nul.
En particulier, dans le contexte galoisien ci-dessus, les noyaux étales 
sauvages $\W_{2i}(L_n)$ attachés aux divers étages de la tour $L_\infty
/L$ ont un $\ell$-rang borné indépendamment de $n$ et de $i$.
\end{Conj}

\Preuve D'après \cite{J1}, les deux groupes $\varprojlim \Cl_{L_n}$ et 
$\C_{F_\infty}$ ont même invariant $\mu$.\medskip

Réunissant les conjectures d'Iwasawa et de Greenberg, 
nous obtenons ainsi~:

\begin{Conj} Soient $F$ un corps de nombres totalement réel puis  $L$ 
comme plus haut abélien sur $F$ et contenant les racines $2\ell$-ièmes 
de l'unité. Alors~:

(i) Les $\ell$-noyaux étales sauvages $\W_{4i}(F_n)$ sont d'ordre borné
indépendamment de $n$ et de $i$ dans la $\Zl$-extension cyclotomique 
$F_\infty = \cup_{n\in\N}\ F_n$.

(ii) Le même résultat vaut pour les composantes imaginaires 
$\W^-_{4i+2}(L_n)$ des $\ell$-noyaux $\W_{4i+2}(L_n)$ dans la $\Zl$-extension cyclotomique 
$L_\infty = \cup_{n\in\N}\ L_n$.
\end{Conj}

\begin{Sco} Les mêmes groupes  $\W_{4i}(F_n)$ (resp. $\W_{4i+2}^-(L_n)$) 
sont ainsi ultime\-ment isomorphes et coïncident avec 
leurs sous-groupes respectifs de capitulation.
\end{Sco}

\Remarque Ce dernier résultat a d'important conséquences heuristiques 
sur les formules de translation du genre {\it à la Riemann-Hurwitz} 
dans les $\ell$-extensions de corps surcirculaires (i.e. de $\Zl$-extensions cyclotomiques 
de corps de nombres)~:\par

Dans un article de synthèse sur ces diverses formules \cite{JMi}, nous 
avons montré, en effet, que leur démonstration se ramène au cas cyclique 
de degré premier, pour lequel  l'information est codée
par le quotient de Herbrand des unités\footnote{Concrètement pour 
le cas des groupes $\C_{N^c}$ dans une extension $\ell$-extension 
cyclique  $N^c\!/F^c$ de groupe $G \simeq \Z/\ell\Z$, par la quantité : $q(G,\wi\E_{N^c})
 = \dim_{\Fl} (H^2(G,\wi\E_{N^c})-\dim_{\Fl} (H^1(G,\wi\E_{N^c})$.}, dont 
l'expérience montre que le calcul  est  impraticable dès qu'il est 
non-trivial. C'est pourquoi, depuis Kida et Iwasawa, on se 
restreint aux seules composantes imaginaires, pour qui le groupe 
des unités se réduit au sous-groupe des racines.\par

Dans \cite{N1}, T. Nguyen Quang Do propose de contourner cette difficulté 
en raisonnant sur les noyaux sauvages. Mais une clef du calcul est 
l'injectivité des morphismes d'extension $\W_{2i}(F_n) \mapsto \W_{2i}
(N_\infty)$, qui suppose bien évidemment la trivialité des noyaux de 
capitulation $C\W_{2i}(F_n)$. Pour les groupes $\W_{4i}(F_n)$, 
l'hypothèse est bien vérifiée d'après le Corollaire 17, mais dans ce 
cas le résultat n'apporte rien de plus que la formule de Kida~; pour les 
groupes $\W_{4i+2}(F_n)$, en revanche, elle ne peut l'être sous les 
conjectures précédentes que si ces groupes sont ultimement triviaux, 
auquel cas le calcul est sans intérêt. En fin de compte, la perspective 
proposée apparait de ce fait quelque peu illusoire eu égard aux conjectures standard.

\smallskip

{\small

\def\refname
}

\bigskip

\begin{tabular}{l p{1.7 cm}  l}
Jean-Fran\c cois {\sc Jaulent}		&{}	& Alexis {\sc Michel}\\
Universit\'e  Bordeaux 	1		&{} 	& Universit\'e  Bordeaux 1\\
Institut de Math\'ematiques		&{} 	& Institut de  Math\'ematiques\\
 351, cours de la Lib\'eration		&{} 	& 351, cours de la Lib\'eration\\
 F-33405 TALENCE  Cedex		&{}	&  F-33405 TALENCE  Cedex\\
{\small jaulent@math.u-bordeaux1.fr}	&{} 	& {\small michel@math.u-bordeaux1.fr}\\
\end{tabular}

\end{document}